\documentclass{amsart}
\usepackage{amsmath}
\usepackage{amsxtra}
\usepackage{amscd}
\usepackage{amsthm}
\usepackage{amsfonts}
\usepackage{amssymb}
\usepackage{eucal}

\newtheorem{cor}{Corollary}
\newtheorem{lem}{Lemma}
\newtheorem{prop}{Proposition}

\newtheorem{thm}{Theorem}

\newcommand{\thmref}[1]{Theorem~\ref{#1}}
\newcommand{\secref}[1]{Sect.~\ref{#1}}
\newcommand{\lemref}[1]{Lemma~\ref{#1}}

\newcommand{\propref}[1]{Proposition~\ref{#1}}

\numberwithin{thm}{section}
\numberwithin{cor}{section}
\numberwithin{prop}{section}
\numberwithin{lem}{section}
\numberwithin{equation}{section}

\newcommand{\nc}{\newcommand}
\nc{\ssec}{\subsection}
\nc{\sssec}{\subsubsection}
\nc{\on}{\operatorname}
\nc{\bi}{\bibitem}

\nc{\ZZ}{{\mathbb Z}}
\nc{\CC}{{\mathbb C}}

\nc{\T}{{\mathcal T}}
\nc{\F}{{\mathcal F}}
\nc{\E}{{\mathcal E}}
\nc{\D}{{\mathcal D}}
\nc{\X}{{\mathcal X}}
\nc{\U}{{\mathcal U}}
\nc{\V}{{\mathcal V}}
\nc{\Y}{{\mathcal Y}}
\nc{\G}{{\mathcal G}}
\nc{\C}{{\mathcal C}}

\nc\Hom{\on{Hom}}

\renewcommand\v{\on{v}}
\nc\w{\on{w}}

\nc\ol{\overline}
\nc\wt{\widetilde}
\nc\tboxtimes{\wt{\boxtimes}}

\nc\grg{\mathfrak g}
\nc\grb{\mathfrak b}
\nc\grt{\mathfrak t}
\nc\grm{\mathfrak m}
\nc\grs{\mathfrak s}

\newcommand\alp{\alpha}		
\newcommand\eps{\varepsilon}		
	
\newcommand\lam{\lambda}
\newcommand\Lam{\Lambda}

\newcommand\calG{{\mathcal{G}}}
\newcommand\calK{{\mathcal{K}}}
\newcommand\calO{{\mathcal{O}}}
\newcommand\calS{{\mathcal{S}}}

\newcommand\bfb{{\mathbf b}}		\newcommand\bfB{{\mathbf B}}
		\newcommand\bfC{{\mathbf C}}
\newcommand\bfd{{\mathbf d}}		
\newcommand\bfe{{\mathbf e}}

\nc\la{\langle}
\nc\ra{\rangle}
\nc\ten{\otimes}
\nc\x{\times}
\newcommand{\gd}{\grg^{\vee}}
\newcommand{\alpv}{\alp^{\vee}}
\renewcommand{\gg}{\calG_G}
\newcommand{\go}{G(\calO)}
\newcommand{\mo}{M(\calO)}
\newcommand{\rc}{\rho^{\vee}}
\newcommand{\ogg}{{\overline \gg}}
\newcommand{\ogl}{\ogg^{\lam}}
\newcommand{\pgg}{\text{Perv}_{\go}(\gg)}
\newcommand{\IC}{{\operatorname{IC}}}
\newcommand{\Res}{\operatorname{Res}}
\newcommand{\Perv}{\operatorname{Perv}}
\newcommand{\rep}{\operatorname{Rep}}
\newcommand{\tilthe}{{\widetilde \theta}}

\begin{document}

\title[Crystals via the the affine grassmannian]{Crystals
via the affine grassmannian}

\author[A.~Braverman and D.~Gaitsgory]{Alexander Braverman and
Dennis Gaitsgory}

\begin{abstract}
Let $G$ be a connected reductive group over $\CC$ 
and let $\gd$ be the Langlands
dual Lie algebra. Crystals for $\gd$ are combinatorial
objects, that were introduced by Kashiwara (cf. for example
\cite{K}) as certain ``combinatorial skeletons'' of finite-dimensional
representations of $\gd$. For every dominant weight $\lam$ of
$\gd$ Kashiwara constructed a crystal $\bfB(\lam)$ by considering the
corresponding finite-dimensional representation of the quantum
group $U_q(\gd)$ and then specializing it to $q=0$. Other
(independent) constructions of $\bfB(\lam)$ were given by 
Lusztig (cf. \cite{Lus2}) using the combinatorics of root
systems and by Littelmann (cf. \cite{Lit}) using 
the ``Littelmann path model''.
It was also shown in \cite{J} that the family
of crystals $\bfB(\lam)$ is unique if certain reasonable conditions
are imposed (cf. \thmref{uniqueness}).

The purpose of this paper is to give another (rather simple) construction
of the crystals $\bfB(\lam)$ using the geometry of the {\it affine
grassmannian} $\calG_G=G(\calK)/G(\calO)$ of the group $G$,
where $\calK=\CC((t))$ is the field of Laurent power series
and $\calO=\CC[[t]]$ is the ring of Taylor series. We then check that
the family $\bfB(\lam)$ satisfies the conditions of the uniqueness 
theorem from \cite{J}, which shows that our crystals coincide with
those constructed in {\it loc. cit}. It would be interesting 
to find these isomorphisms directly (cf., however, \cite{Lus3}).
\end{abstract}
\maketitle

\section{Basic results about crystals}
\ssec{Notation} Let $G$ be a connected reductive
group over $\CC$ and let $G^{\vee}$ be the Langlands dual group;
let $\gd$ denote the Lie algebra of $G^{\vee}$. Let also 
$\text{Rep}(G^{\vee})$ denote the category of 
finite-dimensional representations of the group $G^{\vee}$.

Let $\Lam_G$ denote the {\it coweight} lattice of $G$, which is the
same as the weight lattice of $G^{\vee}$. Let $\Lam_G^{\vee}$ denote
the dual lattice, i.e. $\Lam_G^{\vee}$ is the {\it weight} lattice
of $G$; let $\langle,\rangle$ be the canonical pairing between 
$\Lam_G$ and $\Lam_G^{\vee}$. We will denote by
$\Lam_G^+$ the semi-group of dominant coweights. Let $I$ denote the
set of vertices of the Dynkin diagram corresponding to $G$. 
For $i\in I$ we will denote by $\alp_i\in\Lam_G$ the 
corresponding simple coroot and by $\alpv_i\in\Lam_G^{\vee}$ the 
corresponding simple root. Let $2\rho_G^{\vee}\in \Lam_G^{\vee}$
be the sum of all positive roots of $G$. 
For $\lam_1,\lam_2\in\Lam_G$, we will write
$\lam_1\underset{G}\geq\lam_2$ if $\lam_1-\lam_2$ is a linear combination
of the $\alp_i$ with non-negative coefficients.

Let $E_i,F_i$ (for $i\in I$) denote the Chevalley
generators of $\gd$. For every $\lam\in \Lam_G^+$ we will denote by
$V(\lam)$ the irreducible representation of $\gd$ with highest 
weight $\lam$ and for $\mu\in\Lambda_G$, $V(\lam)_\mu$ will denote the 
corresponding weight subspace of $V(\lam)$.

\ssec{Definition}
A crystal is a set $\bfB$ together with maps
\begin{enumerate}
\item $wt: \bfB\to \Lam_G,\ \eps_i,\phi_i: \bfB\to \ZZ$,
\item $e_i,f_i:\bfB\to \bfB\cup \{ 0\}$,
\end{enumerate}
for each $i\in I$, satisfying the following axioms:

A) For any $\bfb\in \bfB$ one has $\phi_i(\bfb)=\eps_i(\bfb)+\la wt(\bfb),\alpv_i\ra$

B) Let $\bfb\in \bfB$.  If $e_i \cdot \bfb\in \bfB$ for some $i$. Then

$$wt(e_i\cdot \bfb)=wt(\bfb)+\alp_i,\ \eps_i(e_i \cdot \bfb)=\eps_i(\bfb)-1,\ 
\phi_i(e_i\cdot b)=\phi_i(b)+1.$$

If $f_i\cdot \bfb\in \bfB$ for some $i$ then 

$$wt(f_i\cdot \bfb)=wt(\bfb)-\alp_i,\ \eps_i(f_i \cdot \bfb)=\eps_i(\bfb)+1,\ 
\phi_i(f_i\cdot \bfb)=\phi_i(\bfb)-1.$$

C) For all $\bfb,\bfb'\in \bfB$ one has $\bfb'=e_i\cdot \bfb$ if an only if 
$\bfb=f_i\cdot \bfb'$.

\smallskip

\noindent
{\it Remark.} In \cite{J} a more general definition of crystals
is considered, where the maps $\eps_i$ and $\phi_i$ are allowed
to assume infinite values. However, such crystals will never appear in 
this paper.

A crystal is called {\it normal} if one has
\begin{equation} \label{normal}
\eps_i(\bfb)=\max\{ n|\ e_i^n\cdot \bfb\neq 0\},\quad
\phi_i(\bfb)=\max\{n|\ f_i^n\cdot \bfb\neq 0\}
\end{equation}

\smallskip

From now on we will consider only normal crystals. Thus, the maps 
$\eps_i$ and $\phi_i$ will be uniquely recovered from 
$wt$, $e_i$ and $f_i$.

\ssec{Tensor product of crystals}
Let $\bfB_1$ and $\bfB_2$ be two
crystals. Following Kashiwara (\cite{K}) we define their tensor
product $\bfB_1\ten \bfB_2$ as follows. As a set $\bfB_1\ten \bfB_2$ is
just equal to $\bfB_1\x \bfB_2$. The corresponding maps are defined in the
following way. Let $\bfb_1\in \bfB_1, \bfb_2\in \bfB_2$. We will denote 
by $\bfb_1\ten \bfb_2$ be the corresponding element in $\bfB_1\x \bfB_2$. Then we set
$$wt(\bfb_1\ten \bfb_2)=wt(\bfb_1)+wt(\bfb_2),$$
$$e_i\cdot (\bfb_1\ten \bfb_2)=
\begin{cases}
e_i \cdot \bfb_1\ten \bfb_2, \quad \text{if}\quad  \eps_i(\bfb_1)> \phi_i(\bfb_2)\\
\bfb_1\ten e_i \cdot \bfb_2,\quad \text{otherwise}
\end{cases}$$
$$f_i\cdot (\bfb_1\ten \bfb_2)=
\begin{cases}
f_i \cdot \bfb_1\ten \bfb_2, \quad \text{if}\quad  \eps_i(\bfb_1)\geq \phi_i(\bfb_2)\\
\bfb_1\ten f_i \cdot \bfb_2,\quad \text{otherwise}
\end{cases}$$

$$\eps_i(\bfb_1\ten \bfb_2)=
\max\{\eps_i(\bfb_2),\eps_i(\bfb_1)-\phi_i(\bfb_2)+\eps_i(\bfb_2)\}$$

$$\phi_i(\bfb_1\ten \bfb_2)=
\max\{\phi_i(\bfb_1),\phi_i(\bfb_2)-\eps_i(\bfb_1)+\phi_i(\bfb_1)\}.$$

\smallskip

It is known (cf. \cite{J}) that $\bfB_1\ten \bfB_2$ is 
crystal and that $\ten$ is an associative operation on crystals.
Moreover, if $\bfB_1$ and $\bfB_2$ are normal then $\bfB_1\ten \bfB_2$ 
is normal as well.

\ssec{Highest weight crystals}
Let $\bfB$ be a crystal. We say that $\bfB$ is {\it a highest weight crystal
of weight $\lam\in\Lam_G$} if there exists an element $\bfb_\lam\in \bfB$,
such that
\begin{enumerate}
\item $wt(\bfb_\lam)=\lam$.
\item $e_i\cdot \bfb_\lam=0$ for every $i\in I$.
\item $\bfB$ is generated by all the $f_i$ acting on $\bfb_\lam$.
\end{enumerate}

It is clear from \eqref{normal} that if $\bfB$ is a normal crystal,
then one necessarily has $\lam\in \Lam_G^+$. The following lemma
gives a useful reformulation of the definition of a highest weight
crystal.
\begin{lem} \label{highest} 
A  crystal $\bfB$ is a highest weight crystal of highest
weight $\lam$ if and only if there exists an element
$\bfb_\lam\in \bfB$, such that
\begin{enumerate}
\item $wt(\bfb_\lam)=\lam$ and $wt(\bfb)< \lam$ for every $\bfb\in \bfB-\bfb_\lam$.
\item $e_i\cdot \bfb_\lam=0$ for every $i\in I$.
\item For every $\bfb\in \bfB-\bfb_\lam$ there exists $i\in I$
such that $e_i\cdot \bfb\neq 0$.
\end{enumerate}
\end{lem}

\ssec{Closed families of crystals}   \label{closed}
Assume that for every $\lam\in \Lam_G^+$ we are given a normal
crystal $\bfB(\lam)$ of highest weight $\lam$. We say that
the $\bfB(\lam)$ form a closed family of crystals if for every
$\lam,\mu\in\Lam_G^+$ there exists an embedding 
$\bfB(\lam+\mu)\hookrightarrow \bfB(\lam)\ten \bfB(\mu)$ (which necessarily
sends $\bfb_{\lam+\mu}$ to $\bfb_{\lam}\ten \bfb_{\mu}$).

\begin{thm}(cf. \cite{J}, 6.4.21)    \label{uniqueness}
Assume that $G$ is of adjoint type. Then there
exists a unique closed family of crystals $\bfB(\lam)$.
\end{thm}

Different constructions of closed families of crystals were
given by Kashiwara (\cite{K}) using quantum groups and
by Lusztig (\cite{Lus2}) and Littelmann (\cite{Lit}) using the
combinatorics of the root systems. The main goal of this paper
is to give another construction of the closed family $\bfB(\lam)$,
using the geometry of the affine Grassmannian.

\section{Basic results about affine Grassmannian}

\ssec{Definition}      \label{gras-not}
Let $\calK=\CC((t))$, $\calO=\CC[[t]]$. By the {\it affine
Grassmannian} of $G$ we will mean the quotient
$\calG_G=G(\calK)/G(\calO)$. It is known (cf. \cite{BD}) that
$\calG_G$ is the set of $\CC$-points of an ind-scheme over
$\CC$, which we will denote by the same symbol.

The orbits of the group $G(\calO)$ on $\gg$ can be described
as follows. One can identify the lattice $\Lam_G$ with
the quotient $T(\calK)/T(\calO)$. Fix $\lam\in\Lam_G^+$ and
let $\lam(t)$ denote any lift of $\lam$ to $T(\calK)$.
Let $\gg^{\lam}$ denote the $\go$-orbit of $\lam(t)$ 
(which clearly does not depend on the choice of $\lam(t)$).
Then it is well-known (cf. \cite{Lus1}) that
$$\gg=\bigsqcup\limits_{\lam\in\Lam_G^+}\gg^{\lam}.$$
Moreover, for every $\lam\in\Lam_G^+$ the orbit
$\gg^{\lam}$ is finite-dimensional and its dimension is
equal to $\la \lam,2\rc_G\ra$. 

Let $\ogl$ denote the closure of $\gg^{\lam}$ in $\gg$;
this is an irreducible projective algebraic variety.
We will denote by $\IC^{\lambda}$ the intersection
cohomology complex on $\ogl$. Let $\pgg$ denote the category of
$G(\calO)$-equivariant perverse sheaves on $\calG_G$. It is known
that every object of this category is a direct sum of the
$\IC^{\lam}$.

\ssec{The convolution}      \label{convolution}
Define the ind-scheme 
$\gg\star\gg$ to be $G(\calK)\underset{\go}\x\gg$. 
Let $$\pi:G(\calK)\times \gg\to\gg\star\gg$$ denote the natural projection. 
One has the natural maps $p_1,p_2:G(\calK) \times \gg \to\gg$ and 
$m:\gg\star\gg\to\gg$ defined as follows.
Let $g\in G(\calK),x\in\gg$. Then
$$p_1(g,x)=g\, \text{mod}\,\go; \quad p_2(g,x)=x;\quad
m(g,x)=g\cdot x.$$
For $\lam_1,\lam_2\in\Lam_G^+$ we set $\gg^{\lam_1}\star\gg^{\lam_2}=
\pi(p_1^{-1}(\gg^{\lam_1})\cap p_2^{-1}(\gg^{\lam_2}))$. In addition, we define
$$(\gg^{\lam_1}\star\gg^{\lam_2})^{\lam_3}=m^{-1}(\gg^{\lam_3})\cap 
\gg^{\lam_1}\star\gg^{\lam_2}$$
It is known (cf. \cite{Lus1}) that 
\begin{equation} \label{Lusestimate}
\on{dim}((\gg^{\lam_1}\star\gg^{\lam_2})^{\lam_3})=
\la \lam_1+\lam_2+\lam_3,\rc_G\ra.
\end{equation}
(It is easy to see that although $\rc_G\in \frac{1}{2}\Lam_G^{\vee}$,
the RHS of \eqref{Lusestimate} is an integer whenever the above intersection
is non-empty.)

For any $\calS_1,\calS_2\in\pgg$ we define 
the convolution $\calS_1\star\calS_2$ as follows. Consider
$p_1^*\calS_1\ten p_2^*\calS_2$. Then due to the fact that
$\calS_1$ is $\go$-equivariant, there exists a canonical perverse sheaf  
$\calS_1\widetilde{\ten}\calS_2$ on $\gg\star\gg$ such that
$\pi^*(\calS_1\widetilde{\ten}\calS_2)\simeq p_1^*\calS_1\ten p_2^*\calS_2$.

We define 
$$\calS_1\star\calS_2=m_!(\calS_1\widetilde{\ten}\calS_2).$$

\begin{thm}(cf. \cite{Lus1},\cite{Gi} and \cite{MV})       \label{grassmannian}
\begin{enumerate}
\item 
Let $\calS_1,\calS_2\in \Perv_{\go}(\gg)$. Then
$\calS_1\star\calS_2\in \Perv_{\go}(\gg)$.
\item 
The convolution $\star$ extends to a structure of
a tensor category on $\Perv_{\go}(\gg)$, which is
equivalent to the category $\text{Rep}(G^{\vee})$.
\end{enumerate}
\end{thm}

\ssec{Restriction functors to Levi subgroups}  \label{restr}
Let $P$ be a Borel subgroup in $G$ and let $N_P$ be its unipotent radical. 
Let $M=P/N_P$ be the corresponding Levi factor. Let $P^{\vee}$ and
$M^{\vee}$ be the corresponding parabolic and Levi subgroups of $G^{\vee}$.
We have the restriction
functor $\Res^{G^{\vee}}_{M^{\vee}}:\rep(G^{\vee})\to \rep(M^{\vee})$.
Let us explain how to represent this functor geometrically, i.e. 
as a functor $\pgg\to \Perv_{\mo}(\calG_M)$.

Let $\Lam_{G,P}$ denote the lattice of characters of the torus
$Z(M^{\vee})$ (the center of $M^{\vee}$). There is a natural surjection
$\alp_{G,P}:\Lam_G\to\Lam_{G,P}$. One can identify 
$\Lam_{G,P}$ with the set of connected components of 
$\calG_M$.

One can also identify $\Lam_{G,P}$ with the set of
orbits of the group $[P,P](\calK)\cdot M(\calO)$ on 
$\gg$. This is done in the following way.
Let $\theta\in \Lam_{G,P}$. Fix a lift $\tilthe$ of $\theta$ to
$\Lam_G$. Let $S_P^{\theta}$ denote the
$[P,P](\calK)\cdot M(\calO)$-orbit of the element $\tilthe(t)\in T(\calK)$
(cf. \secref{gras-not}). It is easy to see that $S_P^{\theta}$ 
depends only on $\theta$ (and not on the choice of $\tilthe(t)$).

\begin{lem} The following hold:
\begin{enumerate}
\item
One has
$\gg=\bigsqcup\limits_{\theta\in \Lam_{G,P}}S_P^{\theta}$.
\item 
Let $\calG_M^{\theta}$ denote the connected component of $\calG_M$
corresponding to $\theta$. 
Then there exists a canonical $[P,P](\calK)\cdot M(\calO)$--equivariant map
$\grt^{\theta}_P:S_P^{\theta}\to \calG_M^{\theta}$
which is equal to identity on the set
$$
\{\nu\in \Lam_G=T(\calK)/T(\calO)|\ \alp_{G,P}(\nu)=\theta\}.
$$ 
(Note that this set is naturally embedded into both $S_P^{\theta}$ and 
$\calG_M^{\theta}$ due to the fact that $T$ is embedded in both 
$G$ and $M$). 
\end{enumerate}
\end{lem}

Let $\nu\in\Lam_M^+\subset\Lam_G$ and let $\theta=\alp_{G,P}(\nu)$. 
Let us denote by
$S_P^{\nu}$ the pre-image $(\grt_P^{\theta})^{-1}(\calG^\nu_M)
\subset S_P^{\theta}$. The schemes $S_P^{\nu}$ are nothing but orbits of 
the group
$N_P(\calK)\cdot M(\calO)$ on $\gg$. We will denote by $\grt^\nu_P$ the 
restriction of $\grt_P^{\theta}$ to $S_P^{\nu}$.

\begin{thm} (\cite{BD}, cf. also \cite{BG} and \cite{MV})  \label{restriction}
\begin{enumerate}
\item
Let $\nu$ (resp., $\lambda$) be a dominant integral coweight of $M$ 
(resp., of $G$) Then the intersection
$S^\nu_P\cap\calG_G^\lambda$ has dimension 
$\leq \langle\nu+\lambda,\rc_G\rangle$ and hence the fibers 
of the projection
$$
\grt^\nu_P:S^\nu_P\cap\gg^\lambda\to\calG_M^\nu
$$
are of dimension 
$\leq \langle\nu+\lambda,\rc_G\rangle-\langle\nu,2\rc_M\rangle$.
\item
Let $\IC^{\lambda}|_{S_P^\theta}$ denote the $*$-restriction of
$\IC^{\lam}$ to $S_P^{\theta}$. Then
for $\lambda\in\Lambda_G^+$ and $\theta\in\Lambda_{G,P}$, the direct image 
$$
\grt_P^\theta{}_{!}(\IC^{\lambda}|_{S_P^\theta})
[ \langle \theta,2(\rc_G-\rc_M)\rangle]
$$
lives in the cohomological degrees $\leq 0$ (in the perverse t-structure). 
(In the above formula we have used the fact that 
$2(\rc_G-\rc_M)$ naturally belongs to the dual lattice of $\Lambda_{G,P}$.)
\item
The functor $\Perv_{G(\calO)}(\gg)\to \Perv_{M(\calO)}(\calG_M)$ given by
$$\calS 
\mapsto \underset{\theta}\oplus H^0(\grt_P^\theta{}_{!}(\calS|_{S_P^\theta}) 
[\langle \theta,2(\rc_G-\rc_M)\rangle]$$
has a structure of a tensor functor and under 
the equivalence of \thmref{grassmannian} it 
is naturally isomorphic to $\Res^G_M$.
\end{enumerate}
\end{thm}

If $B$ is a Borel subgroup of $G$ then one has $\Lam_G=\Lam_{G,P}$.
In this case for every $\mu\in\Lam_G$ we will write 
$S^{\mu}$ instead of $S_B^{\mu}$. It is clear that for any parabolic $P$, 
$S^\mu$ lies inside $S_P^{\alpha_{G,P}(\mu)}$.

\section{The construction of $\bfB^G(\lam)$}
In this section we will state our two main theorems.
Their proofs will be given in the next two sections.

\ssec{The set $\bfB^G(\lam)$}  \label{intrB}
Let $M$ be as in \secref{restr}. For $\lam\in\Lam_G^+$ and
$\nu\in\Lam_M^+$ we let
$\bfB^G_M(\lam)_\nu$ denote the set of irreducible components of
the intersection $S_P^\nu\cap \calG_G^\lam$
of dimension $\la \nu+\lam,\rc_G\ra$.
Since the variety $\calG_M^\nu$ is connected and simply connected,
it follows, that $\bfB^G_M(\lam)_\nu$ can also be identified
with the set of irreducible components of any fiber of the map
$\grt^\nu_P:S^\nu_P\cap\gg^\lambda\to\calG_M^\nu$ of
dimension  $\la \nu+\lam,\rc_G\ra -\la \nu,2\rc_M\ra$. 

For $\mu\in\Lam_G$ we will denote $\bfB^G_T(\lam)_\mu$ just
by $\bfB^G(\lam)_\mu$ and we set
$$\bfB^G(\lam):=\bigcup\limits_{\mu\in \Lam_G}\bfB^G(\lam)_\mu.$$

Thus, $\bfB^G(\lam)$ is a finite set, endowed with
a map $wt:\bfB^G(\lam)\to \Lam_G$ (by definition, $wt(\bfb)=\mu$ for 
$\bfb\in \bfB^G(\lam)_\mu$).

\ssec{Decomposition with respect to a parabolic}
We would like now to extend the map $wt:\bfB^G(\lam)\to \Lam_G$ to a structure
of a normal crystal on $\bfB^G(\lam)$, i.e. we need to define
the operations $e_i$ and $f_i$.

Let $P$ be any parabolic subgroup in $G$.
\begin{prop}  \label{bijection}
For every $\lam\in \Lam_G^+,\mu\in\Lam_G$ 
there is a canonical bijection
$$\bfd_M^G:  \bigsqcup\limits_{\nu\in\Lam_M^+} \bfB^G_M(\lam)_\nu\x 
\bfB^M(\nu)_\mu\simeq \bfB^G(\lam)_\mu.$$  
This bijection can be uniquely characterized as follows: one has
$\bfd(\bfb_1,\bfb_2)=\bfb$ for $\bfb_1\in \bfB^G_M(\lam)_\nu,
\bfb_2\in \bfB^M(\nu)_\mu$ if  and only if 
the following conditions hold.
\begin{enumerate}
\item $\theta:=\alp_{G,P}(\mu)=\alp_{G,P}(\nu)$.
\item $\bfb_2$ is a dense subset of $\grt_P^{\theta}(\bfb)$.
\item $(\grt_P^{\nu})^{-1}(\bfb_2)\cap \bfb_1$ is a dense subset of 
$\bfb$.
\end{enumerate}
\end{prop}

\begin{proof}

For $\bfb_2\in \bfB^M(\nu)_\mu$ consider the variety 
$(\grt_P^{\nu})^{-1}(\bfb_2)\cap\calG_G^\lam
\subset S^\mu\cap\calG_G^\lam$. It follows from \secref{intrB}
that the set of its irreducible components of
dimension $\la \mu+\lam,\rc_G\ra$ is in a bijection with
$\bfB^G_M(\lam)_\nu$.

Thus, for $\bfb_1\in \bfB^G_M(\lam)_\nu$, we set
$\bfd_M^G(\bfb_1\times \bfb_2)$ to be the closure in $S^\mu\cap\calG_G^\lam$
of the corresponding irreducible component of
$(\grt_P^{\nu})^{-1}(\bfb_2)\cap\calG_G^\lam$.

The fact that this map is a bijection satisfying all the 
required properties is straightforward.

\end{proof}

\ssec{Operations $e_i$ and $f_i$.}
Fix now any $i\in I$. Let $P_i$ be the corresponding
``sub-minimal'' parabolic subgroup of $G$ (by definition,
$P_i$ is the parabolic subgroup of $G$, whose unipotent
radical contains all simple roots except for $\alp_i^{\vee}$).
Let also $M_i$ be the corresponding Levi factor and
$\grm_i^{\vee}$ the dual Lie algebra.

Consider the decomposition of \propref{bijection} for $M=M_i$.
Since $\grm_i^{\vee}$ is a reductive Lie algebra, whose
semi-simple part is isomorphic to ${\bf sl}(2)$, it follows
from \thmref{restriction}(3) and the representation theory
of ${\bf sl}(2)$ that for every $\bfb_2\in \bfB^{M_i}(\nu)_\mu$ there
exists no more than one $\bfb_2'\in \bfB^{M_i}(\nu)_{\mu+\alp_i}$ 
(resp. $\bfb_2''\in \bfB^{M_i}(\nu)_{\mu-\alp_i}$).

Let now $\bfb\in \bfB^G(\lam)_\mu$. Assume that 
$\bfb=\bfd^G_M(\bfb_1\times \bfb_2)$. Thus we define
$$e_i\cdot \bfb=
\begin{cases}
\bfd^G_M(\bfb_1\x \bfb_2')\quad\text{if there exists 
$\bfb_2'\in \bfB^{M_i}(\nu)_{\mu+\alp_i}$}\\
0\quad\text{otherwise}
\end{cases}$$
and
$$f_i\cdot \bfb=
\begin{cases}
\bfd^G_M(\bfb_1\x \bfb_2'')\quad\text{if there exists 
$\bfb_2''\in \bfB^{M_i}(\nu)_{\mu-\alp_i}$}\\
0\quad\text{otherwise}
\end{cases}$$

\begin{thm}        \label{main}
\begin{enumerate}
\item
The maps $e_i$, $f_i$ and $wt$ define a structure of a normal
crystal on $\bfB^G(\lam)$.
\item
The crystal $\bfB^G(\lam)$ defined above is a highest weight
crystal of highest weight $\lam$.
\item
The crystals $\bfB^G(\lam)$ defined above form a closed family (in the sense
of \secref{closed}).
\end{enumerate}
\end{thm}

\medskip

The first point of this theorem follows readily from the representation theory
of ${\bf sl}(2)$. The geometric content of the second point of \thmref{main} 
is summarized in the next corollary:

Let $w_0$ denote the longest element of the Weyl group of $G$ and for 
$i\in I$ let $\grs_i$ be the corresponding
simple reflection. Let $\lam,\mu$ be a pair of elements of $\Lam_G$ with 
$\lam\in\Lam_G^+$. Let
$\bfb$ be an irreducible component of dimension $\la\lam+\mu,
\rc_G\ra$ of $S^\mu\cap\gg^\lam$.

\begin{cor}
Assume that $\mu\neq\lam$ (resp., $w_0(\mu)\neq \lam$). 
Then one can find $i\in I$ and 
$\nu\in\Lam_{M_i}^+$ with $\mu\neq\nu$ 
(resp., $\grs_i(\mu)\neq \nu$) such that the map
$\grt_{P_i}^{\nu}:(\bfb\cap S_{P_i}^\nu)\to S_{M_i}^\mu\cap 
\calG_{M_i}^{\nu}$ is dominant.
\end{cor}

\medskip

Finally, we note that the third point of \thmref{main} combined with
\thmref{uniqueness} implies that our crystals $\bfB^G(\lam)$ are isomorphic 
to those constructed in \cite{K}, \cite{Lus2} and \cite{Lit}.
Indeed, when $G$ is adjoint this is immediate and, in general, if $G$ and $G'$
are isogenous, the corresponding crystals $\bfB^G(\lam)$ and 
$\bfB^{G'}(\lam)$ are isomorphic for $\lam\in\Lam_{G}^+\cap\Lam_{G'}^+$.

\ssec{Refinement}
Here we would like to refine the statement of \thmref{main}(3). Namely,
we want to describe the crystal $\bfB^G(\lam_1)\ten \bfB^G(\lam_2)$ in geometric
terms.

For $\lam_1,\lam_2,\lam_3\in\Lam_G^+$ let 
$\bfC^G(\lam_1,\lam_2)_{\lam_3}$ be the set of all irreducible 
components of dimension $\langle \lam_1+\lam_2+\lam_3,\rc_G\rangle$
of the variety $(\gg^{\lam_1}\star\gg^{\lam_2})^{\lam_3}$.

\begin{thm}   \label{tensor}
One has a canonical isomorphism of crystals
$$\bfB^G(\lam_1)\ten \bfB^G(\lam_2)=\bigsqcup\limits_{\lam_3\in\Lam_G^+}
\bfC^G(\lam_1,\lam_2)_{\lam_3}\x \bfB^G(\lam_3)$$
where the crystal structure on the right hand side comes from
the second multiple.
\end{thm}

\section{Proof of \thmref{main}(2)}

\ssec{Notation}       \label{fundclass}
Let $Z$ be a complex algebraic variety of
dimension $d$ and let $X\subset Z$ be a $d$--dimensional
irreducible component of $Z$. Then we can define
an element $\v(X)\in H^{2d}_c(Z,\CC)$ as follows. Let $Y_1,...,Y_n$
be other irreducible components of $Z$ and let 
$$
X^0=X-\bigcup\limits_{k=1}^n X\cap Y_k.
$$

Denote by $i$ the embedding of $X^0$ into $Z$. Consider the
complex $i_!\CC$ on $Z$. Then, one has a natural map
$i_!\CC\to \CC$ of (complexes of) sheaves on $Z$ and, therefore, a
map
$$H_c^{2d}(Z,i_!\CC)\to H^{2d}_c(Z,\CC).$$

Now, since $X^0$ is irreducible, one has
$$H_c^{2d}(Z,i_!\CC)=H^{2d}_c(X^0,\CC)\simeq \CC.$$

Thus, by composing the above two maps, we get an element
$\v(X)\in H^{2d}_c(Z,\CC)$. Moreover, the collection of elements
$\v(X)$ (for all irreducible components $X$ of $Z$ of the top dimension) 
is a basis of $H^{2d}_c(Z,\CC)$.

\ssec{The basis in $\Hom_M(U(\nu),V(\lam))$}    \label{basis}
Let as before $M$ be a Levi subgroup of $G$. For $\nu\in\Lam_M^+$ we 
will denote by $U(\nu)$ the irreducible representation of
$M$ with highest weight $\nu$. 
We would like now to construct a basis in the vector space
$\Hom_M(U(\nu),V(\lam))$, parametrized by the set $\bfB^G_M(\lam)_\nu$. 
This is done in the following way.

By \thmref{restriction} one can identify $\Hom_M(U_M(\nu),V(\lam))$ with 
$$H_c^{2(\la\lam+\nu,\rc_G\ra-\la\nu,2\rc_M\ra)}
((\grt_P^{\nu})^{-1}(x)\cap \calG_G^\lambda,\CC)$$ 
for any $x\in \calG_M^\nu$.
Recall that $\bfB^G_M(\lam)_\nu$ can be naturally identified with
the set of irreducible components of
$(\grt_P^{\nu})^{-1}(x)\cap \calG_G^\lambda$ of dimension 
$2(\la\lam+\nu,\rc_G\ra-\la\nu,2\rc_M\ra)$.
Hence, the construction of \secref{fundclass} yields a basis
$\v^G_M(\bfb),\,\bfb\in \bfB^G_M(\lam)_\nu$ in $\Hom_M(U(\nu),V(\lam))$.

\ssec{Compatibility of bases}
Fix a weight $\mu\in\Lam_G$ and consider the vector
space $V(\lam)_\mu$. Fix also a parabolic subgroup $P$ with
a Levi subgroup $M$ as before. Then from \secref{basis} one
constructs two bases in $V(\lam)_\mu$, parametrized by $\bfB^G(\lam)_\mu$:
the first one is $\{\v^G_T(\bfb)\}_{\bfb\in B(\lam)_\mu}$ and the other one
is equal to
$$\bigsqcup\limits_{\nu\in\Lam_M^+}\{\v^G_M(\bfb_1)\ten\v^M_T(\bfb_2)|\ 
\text{for }
\bfb_1\in \bfB^G_M(\lam)_\nu\ \text{and }\bfb_2\in \bfB^M(\nu)_\mu\}.$$

Let us now investigate the connection between these two bases.
Let $F^\nu V(\lam)$ denote the direct sum of all $M$-isotypic
components of $V(\lam)$ of the form $U_M(\nu')$, where 
$\nu'\underset{M}\geq \nu$. Set 
$G^\nu V(\lam)=F^\nu V(\lam)/\underset{\nu'\underset{M}>\nu}\sum\,
F^{\nu'} V(\lam)$.

\begin{prop}        \label{triangular}
Let $\bfb\in \bfB^G(\lam)_\mu$. Assume that 
$\bfb=\bfd^G_M(\bfb_1\x \bfb_2)$ where $\bfb_1\in \bfB^G_M(\lam)_\nu$ and
$\bfb_2\in \bfB^M(\nu)_\mu$. Then
\begin{enumerate}
\item $\v^G_T(\bfb)\in F^\nu V(\lam)$.
\item The images of $\v^G_T(\bfb)$ and $\v^G_M(\bfb_1)\ten \v^M_T(\bfb_2)$ 
in $G^\nu V(\lam)$ coincide.
\end{enumerate}
\end{prop}

\begin{proof}

The filtration $F^\nu V(\lam)$ is compatible with the direct sum decomposition
$V(\lam)=\underset{\mu}\oplus V(\lam)_\mu$. Let $F^\nu V(\lam)_\mu$
(resp., $G^\nu V(\lam)_\mu$) denote the corresponding subspace 
(resp., sub-quotient) of $V(\lam)_\mu$.

By \thmref{restriction}, we can identify $V(\lam)_\mu$ with the cohomology
$H_c^{\la \lam+\mu,\rc_G\ra}(S^\mu\cap\gg^{\lam},\CC)$. In addition, 
the filtration $F^\nu V(\lam)_\mu$ on $V(\lam)_\mu$ coincides
with the filtration on the compactly supported cohomology induced by 
the decreasing sequence of open subsets in $S^\mu$:
$$ \underset{\nu'\underset{M}\geq \nu}\bigsqcup 
S^\mu\cap S_P^{\nu'}.$$

Therefore, $G^\nu V(\lam)_\mu\simeq 
H_c^{\la \lam+\mu,\rc_G\ra}(S^\mu\cap S_P^\nu\cap \gg^\lam,\CC)$.

The assertion of the proposition follows now from properties 1--3 of
the bijection $\bfd^G_M$.

\end{proof}

\ssec{Proof of \thmref{main}(2)}
Let us explain how \propref{triangular} implies \thmref{main}(2). 
Conditions 1 and 2 of \lemref{highest} follow from the
well--known fact that the intersection $S^\mu\cap\gg^\lam$ is empty
unless $\lam\geq\mu$ and for $\mu=\lam$, the above intersection 
is dense in $\gg^\lam$ and hence is irreducible. Thus, we just need to prove
that $\bfB^G(\lam)$ satisfies the third condition of \lemref{highest}.

Let $\bfb=\bfd^G_M(\bfb_1\times\bfb_2)$ with $\bfb_1\in\bfB^G_M(\lam)_\nu$ and 
$\bfb_2\in\bfB^M(\nu)_\mu$. Consider the element $\v:=\v^G_T(\bfb)
\in V(\lam)_\mu$. Since $\mu<\lam$,
there exists $i\in I$ and a vector $\v_1\in V(\lam)$ 
such that $F_i(\v_1)=\v$.  
We claim that this implies that $e_i \cdot \bfb\neq 0$.

Indeed, let us denote by $\v'$ the element 
$\v^G_{M_i}(\bfb_1)\ten \v^{M_i}_T(\bfb_2)$. By definition, it is sufficient 
to show that $E_i(\v')\neq 0$. 

We have the canonical $M_i$--invariant projection
$V(\lam)\to G^\nu V(\lam)$ and let $\w$ and $\w'$ be the images 
under this projection of $\v$ and $\v'$, respectively. Now,
\propref{triangular} implies that $\w=\w'$. Hence, if $\w_1$ denotes 
the projection of $\v_1$, we obtain that $F_i(\w_1)=\w'$.
But this means that $E_i(\w')\neq 0$ and hence $E_i(\v')\neq 0$.

\section{Proof of \thmref{tensor}}

\ssec{\thmref{tensor} on the level of sets} 

We will prove a more general assertion. Namely, for a parabolic subgroup $P$
with the Levi factor $M$ and  
$\lam_1,\lam_2\in\Lambda_G^+$ and $\nu\in\Lambda_M^+$,
we will establish a canonical bijection 

\begin{equation}  \label{tensordecomp}
\bfe^G_M:\bigsqcup\limits_{\lam_3\in\Lam_G^+}
\bfC^G(\lam_1,\lam_2)_{\lam_3} \x \bfB^G_M(\lam_3)_\nu\simeq
\bigsqcup\limits_{\nu_1,\nu_2\in\Lam_M^+}
\bfB^G_M(\lam_1)_{\nu_1}\x \bfB^G_M(\lam_2)_{\nu_2}\x
\bfC^M(\nu_1,\nu_2)_\nu
\end{equation}

\begin{proof}

Consider the variety $$m^{-1}(S_P^\nu)\cap (\gg^{\lam_1}\star\gg^{\lam_2}).$$
According to \eqref{Lusestimate} and \thmref{restriction}(1), its
set of irreducible components of dimension
$\la \lam_1+\lam_2+\nu,\rc_G\ra$ can be identified 
with the LHS of \eqref{tensordecomp}.

Now, for $\theta_1,\theta_2\in\Lambda_{G,P}$, let us denote by 
$S_P^{\theta_1}\star S_P^{\theta_2}$ the following scheme:
$$S_P^{\theta_1}\star S_P^{\theta_2}:=[P,P](\calK)M(\calO)\cdot \tilthe_1(t)
\underset{P(\calO)}\times S_P^{\theta_2},$$
where $\tilthe_1(t)$ is as in \secref{restr}. It is easy to see that
the natural map $S_P^{\theta_1}\star S_P^{\theta_2}\to \gg\star\gg$
is a locally closed embedding.

Similarly, for $\nu_1,\nu_2\in\Lambda^+_M$ and 
$\lam_1,\lam_2\in\Lambda^+_G$ we define the sub-scheme
$(S_P^{\nu_1}\cap\gg^{\lam_1})\star (S_P^{\nu_2}\cap\gg^{\lam_2})$ of
$\gg\star\gg$ as 
$S_P^{\theta_1}\star S_P^{\theta_2}\cap \gg^{\lam_1}\star\gg^{\lam_2}$.

We have a commutative diagram
$$
\CD
S_P^{\theta_1}\star S_P^{\theta_2}  @>m>> S_P^{\theta_1+\theta_2}  \\
@V{\grt_P^{\theta_1}\star \grt_P^{\theta_2}}VV   
@V{\grt_P^{\theta_1+\theta_2}}VV  \\
\calG_M^{\theta_1}\star \calG_M^{\theta_2}  @>m>>
\calG_M^{\theta_1+\theta_2}.
\endCD
$$

Therefore, to each element of the set 
$\bfB^G_M(\lam_1)_{\nu_1}\x \bfB^G_M(\lam_2)_{\nu_2}\x
\bfC^M(\nu_1,\nu_2)_\nu$ we can attach an irreducible component of dimension
$\la \lam_1+\lam_2+\nu,\rc_G\ra$ in
$(S_P^{\nu_1}\cap\gg^{\lam_1})\star (S_P^{\nu_2}\cap\gg^{\lam_2})$.
By taking its closure in 
$m^{-1}(S_P^\nu)\cap (\gg^{\lam_1}\star\gg^{\lam_2})$
we obtain an irreducible component of 
$m^{-1}(S_P^\nu)\cap (\gg^{\lam_1}\star \gg^{\lam_2})$
and it is easy to see that the map we have just described is a bijection.

This proves our assertion.

\end{proof}

Note now that for the torus $T$, $\bfC^T(\mu_1,\mu_2)_\mu=\emptyset$ unless
$\mu_1+\mu_2=\mu$ and in the latter case this is the set of one element.
Therefore, for $M=T$ \eqref{tensordecomp} yields the needed isomorphism
$$\bfB^G(\lam_1)\times \bfB^G(\lam_2)\overset{\bfe^G_T}\simeq 
\bigsqcup\limits_{\lam_3\in\Lam_G^+}
\bfC^G(\lam_1,\lam_2)_{\lam_3} \x \bfB^G(\lam_3).$$

\ssec{Compatibility of decompositions}   \label{compatibility}

Consider the set $\bfB^G(\lam_1)\times \bfB^G(\lam_2)$
which, as we have seen above, can be canonically identified with 
$\bigsqcup\limits_{\lam_3\in\Lam_G^+}
\bfC^G(\lam_1,\lam_2)_{\lam_3} \x \bfB^G(\lam_3)$.

There are {\it a priori} two different ways to identify this set with 
$$\bigsqcup\limits_{\nu_1,\nu_2\in\Lambda_M^+} 
\bfB^G_M(\lam_1)_{\nu_1}\x \bfB^G_M(\lam_2)_{\nu_2}\x \bfB^M(\nu_1) 
\x \bfB^M(\nu_2):$$

One is 
$$\bfB^G(\lam_1)\times \bfB^G(\lam_2)\overset{\bfd^G_M\x \bfd^G_M}\longrightarrow
\bigsqcup\limits_{\nu_1,\nu_2\in\Lambda_M^+}
\bfB^G_M(\lam_1)_{\nu_1}\x \bfB^G_M(\lam_2)_{\nu_2}\x 
\bfB^M(\nu_1)\x \bfB^M(\nu_2).$$

The other one is the composition
\begin{align*}
&\bigsqcup\limits_{\lam_3\in\Lam_G^+}
\bfC^G(\lam_1,\lam_2)_{\lam_3} \x \bfB^G(\lam_3) \overset{\bfd^M_T}\simeq
\bigsqcup\limits_{\lam_3\in\Lam_G^+;\nu\in\Lam^+_M}
\bfC^G(\lam_1,\lam_2)_{\lam_3}\x \bfB^G_M(\lam_3)_\nu\times \bfB^M(\nu) 
\overset{\bfe^G_M}\simeq \\
&\bigsqcup\limits_{\nu_1,\nu_2\in\Lambda_M^+}
\bfB^G_M(\lam_1)_{\nu_1}\x \bfB^G_M(\lam_2)_{\nu_2}\x 
\bfC^M(\nu_1,\nu_2)_\nu\x \bfB^M(\nu)\overset{\bfe^M_T}\simeq \\
&\bigsqcup\limits_{\nu_1,\nu_2\in\Lambda_M^+}
\bfB^G_M(\lam_1)_{\nu_1}\x \bfB^G_M(\lam_2)_{\nu_2}\x
\bfB^M(\nu_1)\times \bfB^M(\nu_1).
\end{align*}

However, it is easy to see from the construction
that these two identifications coincide.

\ssec{Reduction to $PGL(2)$}

We have established the isomorphism of sets
$$\bigsqcup\limits_{\lam_3\in\Lam_G^+}
\bfC^G(\lam_1,\lam_2)_{\lam_3} \x \bfB^G(\lam_3)\simeq
\bfB^G(\lam_1)\times \bfB^G(\lam_2)$$
and we must show that the $e_i$ and $f_i$ operations on both sides coincide.

For $i\in I$ consider the corresponding parabolic $P_i$. We decompose
the LHS as 
$$\bigsqcup\limits_{\nu_1,\nu_2\in\Lam^+_{M_i}}
(\bfB^G_M(\lam_1)_{\nu_1}\x \bfB^G_M(\lam_2)_{\nu_2})\x 
(\bfC^M(\nu_1,\nu_2)_\nu\x \bfB^M(\nu))$$
and the RHS as
$$\bigsqcup\limits_{\nu_1,\nu_2\in\Lam^+_{M_i}}
(\bfB^G_M(\lam_1)_{\nu_1}\x \bfB^G_M(\lam_2)_{\nu_2})\x
(\bfB^M(\nu_1)\x \bfB^M(\nu_2)).$$

According to \secref{compatibility}, these decompositions are compatible.
By definition, in both cases, the $e_i$ and $f_i$ operations
preserve these decompositions and act ``along'' the second multiple.

This observation reduces the assertion of \thmref{tensor} from $G$ to
$M_i$. In addition, it is easy to see that we can replace $M_i$ by
its adjoint group, i.e. it remains to analyze the case of $G=PGL(2)$.

\ssec{Proof of \thmref{tensor} for $PGL(2)$}

For $G=PGL(2)$ we will identify $\Lambda_G$ (resp., $\Lam_G^+$)
with $\ZZ$ (resp., with $\ZZ^+$). The positive root $\alpha\in\Lambda_G$
corresponds to $2\in\ZZ$.

Let $l_1,l_2$ be two elements of $\ZZ^+$. The action of $e$ and $f$ breaks 
$\bfB^G(l_1)\otimes \bfB^G(l_2)$ into orbits and it is sufficient to show that
this decomposition coincides with 
$$\bfB^G(l_1)\otimes \bfB^G(l_2)\simeq \bigsqcup\limits_{l\in\ZZ^+} \bfC^G(l_1,l_2)_l\times \bfB^G(l)$$
(note that in this case each $\bfC^G(l_1,l_2)_l$ has at most one element.)

For that end, it is sufficient to show that for
$m_1,m_2\in\ZZ$ a generic point
in $(S^{m_1}\cap \gg^{l_1})\star (S^{m_2}\cap \gg^{l_2})$
projects under the map $m:\gg\star\gg\to\gg$ to $S^n$, where
$$n=\max\{l_1-m_2,m_1+l_2\}.$$

For $l\in\ZZ^+$ and $m\in\ZZ$, the intersection $S^m\cap \gg^l$
is non--empty if only if $l-m\in 2\ZZ^+,\,l\geq |m|$ and in the latter case it
consists of cosets of the form
$$
\begin{pmatrix} t^m & t^{(m-l)/2}\cdot p(t) \\ 0 & 1
\end{pmatrix} \cdot PGL(2,\calO) \mid p(t)\in\CC[[t]],\, p(0)\neq 0 .
$$

Therefore, the image of $(S^{m_1}\cap \gg^{l_1})\star (S^{m_2}\cap \gg^{l_2})$
under $m$ consists of cosets of the form
$$
\begin{pmatrix} t^{m_1+m_2} & t^{\max\{l_1-m_2,m_1+l_2\}}\cdot p(t)
\\ 0 & 1
\end{pmatrix} \cdot PGL(2,\calO) \mid p(t)\in\CC[[t]],\, p(0)\neq 0.
$$

This finishes the proof of \thmref{tensor}. 

\bigskip

\noindent{\bf Acknowledgements.}

We wish to use this occasion in order to thank A.~Joseph for 
explaining to us the basics of crystals. In addition, D.G. would like to thank
R.~MacPherson and J.~Anderson, who were interested in the same problem
from a slightly different angle, for an illuminating discussion.

\end{document}